\documentclass[12pt]{article}
\usepackage{amsfonts}
\usepackage{amssymb}
\usepackage{amsmath}
\usepackage{amsthm}
\usepackage{verbatim}

\newcommand{\delo}{\partial \Omega}
\renewcommand{\div}{\mathrm{div}}

\newcommand{\N}{{\mathbb N}}
\newcommand{\R}{{\mathbb R}}
\newcommand{\Z}{{\mathbb Z}}

\normalbaselines 

\newtheorem{theorem}{Theorem}[section]
\newtheorem{lemma}[theorem]{Lemma}

\newtheorem{definition}{Definition}[section]

\theoremstyle{definition}

\begin{document}

\title{BMO solvability and the $A_\infty$ condition for elliptic operators}

\author{Martin Dindos\footnote{supported by EPRC grant EP/F014589/1-253000} \and Carlos Kenig\footnote{supported by NSF} \and Jill Pipher\footnote{supported by NSF}}

\maketitle

\begin{abstract}
We establish a connection between the absolute continuity of
elliptic measure associated to a second order divergence form
operator with bounded measurable coefficients with the solvability
of an endpoint $BMO$ Dirichlet problem. We show that these two
notions are equivalent. As a consequence we obtain an end-point
perturbation result, i.e., the solvability of the $BMO$ Dirichlet
problem implies $L^p$ solvability for all $p>p_0$.
\end{abstract}

\section{Introduction}

We shall prove an equivalence between
solvability of certain endpoint ($BMO$) Dirichlet boundary value
problems for second order elliptic operators and a quantifiable absolute
continuity of the elliptic measure associated to these operators.
More precisely, we consider here the Dirichlet problem for divergence form (not necessarily
symmetric) elliptic operators $L = \text{div} A \nabla$, where
$A = (a_{ij}(X))$ is a matrix of bounded measurable functions for which
there exists a $\lambda >0$ such that
$\lambda^{-1}\vert \xi \vert^2 < \sum a_{ij}\xi_i\xi_j < \lambda \vert \xi \vert^2$.
The $L^p$ Dirichlet problem for $L$ asks for solvability in a domain
$\Omega$, in the sense
of non-tangential convergence and a priori $L^p$ estimates, of the
problem:
$Lu = 0$ in $\Omega$ with $u=f$ on $\partial \Omega$.

Let us recall (\cite{FS}) a fundamental property of the
harmonic extension to $\R_n^+$ of functions of bounded mean oscillation
on $\R^n$: If $f \in BMO$, then the Poisson extension $u(x,t) = P_t * f(x)$
has the property that $t\vert \nabla u \vert^2 dxdt$ is a Carleson measure.
(Carleson measures are defined in Section 2, below.) In fact the Carleson measure norm of this
extension and the $BMO$ norm of $f$ are equivalent.

In \cite{FN}, this fundamental property was shown to hold for the
harmonic functions in the class of Lipschitz domains.
The key fact here is that harmonic measure on Lipschitz
domains is always mutually absolutely continuous with respect to surface measure, by
a well known result of \cite{D}.

%To establish the context for the results here, let us begin
%by recalling the methods for establishing absolute continuity
%which appeared in \cite{KKPT}. There, the following was shown.
%Let $L = \div A \nabla$ be an elliptic operator with bounded measurable
%coefficients whose matrix is real but not necessarily symmetric.
%If all bounded solutions possess a certain approximation
%property in a domain $D$, then the elliptic measure is in fact absolutely
%continuous in $D$. The approximation property is one which is
%well known for harmonic functions in the disc or the upper half
%space.

In \cite{KKPT}, further connections between Carleson measure
properties of solutions to very general second order divergence
form elliptic equations and absolute continuity were established.
There it was shown that if all bounded solutions
to $L = \text{div} A \nabla$ are arbitrarily well
approximated by continuous functions satisfying an
$L^1$ version of the Carleson measure property, then
in fact the elliptic measure belongs to $A_{\infty}$ with
respect to surface measure.
This approximation property was shown (in \cite{KKPT}) to
follow from a certain norm equivalence between two different classical
quantities associated to the solution of an elliptic equation:
the nontangential maximal function, measuring size, and
the square function, measuring the size of oscillations.

These results, from the Carleson measure properties of
harmonic functions in the upper half space, to theorems such as those
in \cite{KKPT} which specifically connect absolute continuity of the
representing measures associated to
second order divergence form operators to Carleson measure
conditions, led us to a conjecture
concerning solvability of the Dirichlet problem with data in $BMO$.

Specifically, we are interested in properties of the
elliptic measure of an operator $L = \text{div} A \nabla$
which determine that it belongs to the Muckenhoupt
$A_{\infty}$ class with respect to the surface measure on the boundary of
the domain of solvability. On the one hand, $A_{\infty}$ is a ``perturbable" condition, in the
sense that $A_{\infty} = \bigcup A_p = \bigcup B_p$. And when the density of
harmonic measure with respect to surface measure belongs to $B_p$, it turns
out that the Dirichlet problem is solvable with data in $L^q$, where
$1/q + 1/p = 1$. (Again, see section 2 for the definitions.) On the other hand,
a boundary value problem which is equivalent to $A_{\infty}$ would have
to be ``perturbable" as well: solving it would have to imply solvability of
the Dirichlet problem in some $L^q$.
Clearly $L^{\infty}$ cannot be such a perturbable endpoint
space: all solutions satisfy a maximum
principle, a precise version of the $L^{\infty}$ Dirichlet problem.
In the end, perturbing from a $BMO$ problem seems quite natural.

We will use a variety of properties of solutions to divergence
form elliptic operators with bounded measurable coefficients.
The De Giorgi-Nash-Moser theory of the late 1950s and early 1960s
assures us that
weak solutions to these equations are in fact Holder continuous.
Further properties of solutions, of the elliptic measure whose existence
is guaranteed by the maximum principle and the Riesz representation theorem,
and of the relationship of this measure to the Green's function were developed
in the 1970s and 1980s.
For the
basic properties of solutions to divergence form
operators with bounded measurable coefficients, as in
\cite{N} or \cite{CFMS}, one can
consult the introduction of \cite{KKPT} where many
primary references are cited, and where the issues
for the non-symmetric situation are discussed.

\section{Definitions and Statements of Main Theorems}

Let us begin by introducing Carleson measures and square functions
on domains which are locally given by the graph of a function. We
shall assume that our domains are Lipschitz, even though it is
possible to formulate and prove these results with less stringent
geometric conditions on the domain. Most likely, the minimal
geometric conditions required would be chord-arc and
nontangentially accessible. \footnote{This was pointed out to us
by M. Badger.}

\begin{definition}\label{def1}
$\Z \subset \R^n$ is an $M$-cylinder of diameter $d$ if there
exists a coordinate system $(x,t)$ such that
\[
\Z = \{ (x,t)\; : \; |x|\leq d, \; -2Md \leq t \leq 2Md \}
\]
and for $s>0$,
\[
s\Z=\{(x,t)\;:\; |x|<sd, -2Md \leq t \leq 2Md \}.
\]
\end{definition}

\begin{definition}
$\Omega\subset \R^n$ is a Lipschitz domain with Lipschitz
`character' $(M,N,C_0)$ if there exists a positive scale $r_0$ and
at most $N$ cylinders $\{Z_j\}_{j=1}^N$ of diameter $d$, with
$\frac{r_0}{C_0}\leq d \leq C_0 r_0$ such
that\\
\noindent (i) $8Z_j \cap \delo$ is the graph of a Lipschitz function $\phi_j$,
\[
\| \phi_j \|_\infty \leq M, \; \phi_j(0)=0,
\]
\noindent (ii)
\[
\delo=\bigcup_j (Z_j \cap \delo )
\]
\noindent (iii)
\[
Z_j \cap \Omega \supset \left\{ (x,t) \; : \; |x|<d, \;
\mathrm{dist}\left( (x,t),\delo \right) \leq \frac{d}{2}\right\}.
\]

If $Q\in \delo$ and
\[
B_r(Q)= \{ x: |x-Q| \leq r \}
\]
then $\Delta_r(Q)$ denotes the surface ball $B_r(Q)\cap \delo$
and $T(\Delta_r)=\Omega\cap B_r(Q)$ is the called the
Carleson region above
$\Delta_r(Q)$.
\end{definition}

\begin{definition}
Let $T(\Delta_r)$ be a Carleson region associated to a surface
ball $\Delta_r$ in $\delo$. A measure $\mu$ in $\Omega$ is
Carleson if there exists a constant $C=C(r_0)$ such that for all
$r\le r_0$,
\[\mu(T(\Delta_r))\le C \sigma (\Delta_r).\]

For such measure $\mu$ we denote by $\|\mu\|_{Car}$ the number
$$\|\mu\|_{Car}=\sup_{\Delta\subset\delo}\left(\sigma(\Delta)^{-1}\mu(T(\Delta))\right)^{1/2}.$$
\end{definition}

\begin{definition}
A cone of aperture $a$ is a non-tangential approach region for $Q \in \partial \Omega$ of the form
\[\Gamma(Q)=\{X\in \Omega: |X-Q|\le a\;\; dist(X,\partial \Omega)\}.\]
Sometimes it is necessary to truncate the height of $\Gamma$ by
$h$. Then $\Gamma_{h}(Q)= \Gamma(Q) \cap B_h(Q)$.
\end{definition}

We remind the reader that $L$ will stand for $L = \div A \nabla$
where the matrix $A$ has bounded measurable coefficients $a_{i,j}$
and is strongly elliptic: there exists $\lambda$ such that for all
$\xi\in\R^n\backslash\{0\}$,
$$
\lambda |\xi|^2\le \sum a_{i,j}\xi_i\xi_j \le \lambda^{-1}|\xi|^2.
$$

\begin{definition}
If $\Omega \subset \R^n$, and $u$ is a solution to
$L$,  the square function in $Q \in \partial
\Omega$ relative to a family of cones $\Gamma$ is
\[Su(Q)=\left( \int_{\Gamma(Q)} |\nabla u (X)|^2
\delta(X)^{2-n} dX \right)^{1/2}.\] and the non-tangential maximal
function at $Q$ relative to $\Gamma$ is
\[Nu(Q)=\sup \{|u(X)|:X\in \Gamma(Q)\}.\]
Here $\delta(X)=\mbox{dist}(X,\delo)$. We also consider truncated
versions of these operators which we denote by $S_hu(Q)$ and
$N_h(Q)$, respectively; the only difference in the definition is
that the nontangential cone $\Gamma(Q)$ is replaced by the
truncated cone $\Gamma_h(Q)$.
\end{definition}

\begin{definition}
\label{Dir-p}
The Dirichlet problem with the $L^p(\delo,d\sigma)$ data is
solvable for $L$ if the solution $u$ for continuous boundary data
$f$ satisfies the estimate
\begin{equation}
\|N(u)\|_{L^p(\partial \Omega, d\sigma)} \lesssim
\|f\|_{L^p(\partial \Omega, d\sigma)},\label{dp}
\end{equation}
where the implied constant does not depend on the given function.
\end{definition}

\begin{definition}
If $d\mu$ and $d\nu$ are finite measures on the boundary of
$\Omega$, then $d\mu$ belongs to $A_{\infty}$ with respect to
$d\nu$ if for all $\epsilon$ there exists an $\eta$ such that, for
every surface ball $\Delta$ and subset $E \subset \Delta$,
whenever $\nu(E)/\nu(\Delta) < \eta$, then $\mu(E)/\mu(\Delta) <
\epsilon$.
\end{definition}

This space was investigated in \cite{CFe}, where various equivalent definitions were
given. In particular, $d\mu \in A_{\infty}(d\nu)$ if and only if
$d\nu \in A_{\infty}(d\mu)$.

Let us specialize this definition to the domain $\Omega$, to
surface measure $d\sigma$ and to the elliptic measure $d\omega_L$
associated to some divergence form operator $L$. We are assuming
that $d\omega_L$ is evaluated at some fixed point $P$ in the
interior of $\Omega$ so that a solution to $L$ with continuous
data $f$ at the point $P$ is represented by this measure: this
means that $u(P) = \int_{\partial \Omega} f(y) d\omega_L(y)$. If
$d\omega$ belongs to $A_{\infty}(d\sigma)$, then there is a
density function: $d\omega_L(y) = k(y)d\sigma$. The apriori
estimate of definition \ref{Dir-p} turns out to be equivalent to
the fact that the density $k(y)$ satisfies a reverse H\"older
estimate $B_{p'}$. For general $q > 1$, the density $k$ is said to
belong to $B_q(d\sigma)$ if there exists a constant $C$ such that
for every surface ball $\Delta$, $((\sigma(\Delta))^{-1}
\int_{\Delta} k^q d\sigma)^{1/q} < C \sigma(\Delta))^{-1}
\int_{\Delta} k d\sigma$. The relationship between the reverse
H\"older classes and $A_{\infty}$ is (\cite{G} and \cite{CFe})
$$A_\infty(\partial\Omega,d\sigma)=\bigcup_{p>1}B_q(\partial\Omega,d\sigma),$$

\begin{definition} We say that a function $f:\delo\to \R$ belongs
to BMO with respect to the surface measure $d\sigma$ if
$$\sup_{I\subset \delo}\sigma(I)^{-1}\int_I |f-f_I|^2
d\sigma<\infty.$$ Here $f_I=\sigma(I)^{-1}\int_{I}fd\sigma$. We
denote by $\|f\|_{BMO(p)}$ the number
$$\|f\|_{BMO(p)}=\sup_{I\subset \delo}\left(\sigma(I)^{-1}\int_I
|f-f_I|^p d\sigma\right)^{1/p}.$$ It can be shown for any $1\le
p<\infty$ that $\|f\|_{BMO(2)}<\infty$ if and only if
$\|f\|_{BMO(p)}<\infty$. Moreover, $\|.\|_{BMO(p)}$ and
$\|.\|_{BMO(2)}$ are equivalent in the sense that there is a
constant $C>0$ such that the inequality
\begin{equation}
C^{-1}\|f\|_{BMO(p)}\le \|f\|_{BMO(2)}\le
C\|f\|_{BMO(p)}\label{88}
\end{equation}
holds for any $BMO$ function $f$.
\end{definition}

This definition can be modified further. Instead of using the
difference $f-f_I$ in the definition of the BMO norm one can take
\begin{equation}
\|f\|_{BMO^*(p)}=\sup_{I\subset \delo}\inf_{c_I\in
\R}\left(\sigma(I)^{-1}\int_I |f-c_I|^p
d\sigma\right)^{1/p}.\label{8}
\end{equation}
Again, it can be shown that this gives an equivalent norm, i.e.,
there is $C>0$ such that
$$C^{-1}\|f\|_{BMO^*(p)}\le \|f\|_{BMO(2)}\le C\|f\|_{BMO^*(p)}.$$

\begin{definition}\label{d-bmo}
The BMO-Dirichlet problem  is solvable for $L$ if the solution $u$
for continuous boundary data $f$ satisfies
$$\||\nabla u|^2\delta(X)dX\|_{Car}\lesssim \|f\|_{BMO(2)}.$$
Equivalently, there exists a constant $C$ such that for all
continuous $f$,
\begin{equation}
\sup_{\Delta\subset\delo}\sigma(\Delta)^{-1}\iint_{T(\Delta)}|\nabla
u|^2\delta(X)dX \le C \sup_{I\subset \delo}\sigma(I)^{-1}\int_I
|f-f_I|^2 d\sigma. \label{equiv}
\end{equation}
\end{definition}

\noindent{\bf Remark 2.1.} It follows from our results that even
though we define BMO-solvability in the Definition \ref{d-bmo}
only for {\it continuous} boundary data, the solution can be
defined for any BMO function $f:\partial\Omega\to\R$ and moreover
the estimate (\ref{equiv}) will hold. In addition, such a solution
$u$ will have a well-defined nontangential maximal function $N(u)$
for almost every point $Q\in\partial\Omega$ and in the
nontangential sense
$$f(Q)=\lim_{X\to Q,\, X\in\Gamma(Q)} u(X),\qquad\mbox{for a.e. }Q\in\partial\Omega.$$

\bigskip

We now state our main results.

\begin{theorem}\label{theorem:mainth} Let $\Omega$ be a Lipschitz domain and $L$ be a divergence form
elliptic operator with bounded coefficients satisfying the strong
ellipticity hypothesis.

If the elliptic measure $d\omega_L$ associated with $L$ is in
$A_\infty(\delo,d\sigma)$ then the BMO-Dirichlet problem is
solvable for $L$, with in fact equivalence of the two norms
in the estimate (\ref{equiv}).

Conversely, if the estimate (\ref{equiv}) holds for all continuous
functions $f$ with constants only depending on the Lipschitz
character of the domain $\Omega$ and the ellipticity constant of
$L$, then the elliptic measure $d\omega_L$ belongs to
$A_\infty(\delo,d\sigma)$.
\end{theorem}

\noindent{\bf Remark 2.2.} The closure of continuous functions in
BMO norm is the VMO class (\cite{Sa}). From the proof of the
theorem, we will see that $A_\infty$ is actually equivalent to
solvability of a VMO-Dirichlet problem. \vglue3mm

Recall that if a Dirichlet problem for an elliptic operator $L$ is
$L^p$ solvable for some $p\in (1,\infty)$, then it is solvable for
all $L^q$ $p-\varepsilon<q<\infty$, which shows that the
\lq\lq solvability" is stable under small perturbations.

Theorem \ref{theorem:mainth} implies the same kind of stability result for
the end-point BMO problem on the $L^p$ interpolation scale.

\begin{theorem}\label{theorem:mainth2} (Stability of BMO solvability) Let $\Omega$ be a Lipschitz domain and $L$ be a divergence form
elliptic operator with bounded coefficients satisfying the strong
ellipticity hypothesis.

Assume that the estimate (\ref{equiv}) holds. Then there exist
$p_0>1$ such that the $L^p$ Dirichlet problem for $L$ is solvable
for all $p_0<p<\infty$.
\end{theorem}

%%%%%%%%%%%%%%%%%%%%%%%%%%%%%%%%%%%%%%%%%%%%%%%%%%%%%%%%%%%%%%%%%%%%%%%%%%%%%%%%%%%%%%%%%%%%%%%%%%%%%%%%

\section{Proofs}

\noindent{\it We start by proving Theorem
\ref{theorem:mainth}.}\vglue2mm

\noindent We establish the $A_{\infty}$ property of $d\omega_L$ by
assuming the estimate (\ref{equiv}) holds uniformly for continuous
data.

The elliptic measure for $L$ will be abbreviated $d\omega$ and is
evaluated at a fixed interior point, $P_0$, of the domain
$\Omega$.

Let $\Delta$ be a surface ball on the boundary of $\Omega$ of
radius $r$. Let $\Delta'$ be another surface ball of radius $r$
separated from $\Delta$ by a distance of $r$. By assumption, if
$Lu=0$ and $u = f$ on the boundary, we have

\begin{equation}
\label{bmoestimate} \sigma(\Delta')^{-1}\iint_{T(\Delta')}|\nabla
u|^2\delta(X)dX \lesssim \| f \|_{BMO}.
\end{equation}

Let us now assume that $f$ is a positive and continuous function supported in
$\Delta$.

Recall that $S_ru(Q)$ denotes the square function defined using cones
truncated at height $r$. We claim that there exists a constant
$C$ such that for all $Q \in \Delta'$,

\begin{equation}
\label{claim} \omega(\Delta)^{-1}\int_{\Delta} f d\omega \leq C
S_ru(Q)
\end{equation}

To establish this claim, we introduce a little more
notation.

For $Q \in \Delta'$,
set $\Gamma_j(Q) = \Gamma(Q) \cap B_{2^{-j}r}(Q)
\backslash B_{2^{-j-1}r}(Q)$, a slice of the cone $\Gamma(Q)$ at height
$2^{-j}r$.

By Lemma 5.8 (see also 5.13) of \cite{KPi}, we have the following
Poincare type estimate, which was established using Sobolev
embedding and boundary Cacciopoli to exploit the fact that $u$
vanishes on $\Delta'$:

\begin{equation}
(2^jr)^{-2} \int_{\Gamma_j(Q)} u^2 dX \lesssim
\iint_{\Gamma_j(Q)}|\nabla u|^2\delta(X)dX
\end{equation}

Let $A'$ denote a point in $T(\Delta')$ whose distance to the
boundary of $\Omega$ is approximately $r$. By the comparison
theorem for solutions which vanish at the boundary, and with
$G(X)$ denoting the Green's function for $L$ with pole at $P_0$
in $\Omega$,

\begin{equation}
\frac{u(X)}{G(X)} \approx \frac{u(A')}{G(A')},
\end{equation}

for all $X \in \Gamma(Q) \cap T(\Delta')$.

We use this to estimate the square function:

\begin{eqnarray}
S_r^2u(Q) &\ge& \sum_{j=0}^{\infty}
\int_{\Gamma_j}\delta(X)^{2-n}|\nabla u|^2 dX
\\&\ge& \frac{u^2(A')}{G^2(A')} \sum_j^{\infty} (2^{-j}r)^n
\int_{\Gamma_j} G^2(X) dX.
\end{eqnarray}

Now let $A_j$ be a nontangential point in $\Gamma_j$, so that
$|A_j - Q| \approx 2^{-j}r$. By Harnack, $G(X) \approx G(A_j)$ for
all $X \in \Gamma_j(Q)$. Moreover, again by Harnack, there is
constant $C > 1$ for which $G(A_{j-1}) < CG(A_j)$. Thus,

\begin{equation}
\frac{u^2(A')}{G^2(A')} \sum_{j=0}^{\infty} G^2(A_j) \lesssim  S_r^2u(Q),
\end{equation}

and now since $G(A') \leq C^j G(A_j)$, we can sum this series
and we find that
\begin{equation}
u^2(A') \lesssim S^2_ru(Q).
\end{equation}

Since, by properties of harmonic measure, we also know
that $u(A') \approx \omega(\Delta)^{-1}\int_{\Delta} f d\omega$, this
proves \ref{claim}.

\smallskip

For any such $f$, positive, continuous and supported in $\Delta$, the estimate in
\ref{bmoestimate} implies that, for some constant $C_0$,

\begin{equation}
\label{eqn:normbound} (\omega(\Delta)^{-1}\int_{\Delta} f
d\omega)^2 \leq C_0^2 \|f\|_{BMO}^2.
\end{equation}

We now establish absolute continuity of the elliptic measure.
Suppose that $\sigma(\Delta) = r$ and that $\epsilon$ is given.
Let $E \subset \Delta$ be an open set. We shall
find an $\eta$ such that $\sigma(E)/\sigma(\Delta) < \eta$ implies that
$\omega(E)/\omega(\Delta) < \epsilon$.

Let $h = \chi_E$, the characteristic function of $E$.
If $M(h)$ is the Hardy-Littlewood maximal function of $h$ with respect to
surface measure on the boundary of $\Omega$, define (as in \cite{JJ})
the $BMO$ function

\begin{equation}\label{eq:specialbmofunction}
f = \text{max}\{0, 1+\delta\log M(h)\},
\end{equation}

\noindent where $\delta$ is to be determined. The function $f$ has
a structure which is typical of $BMO$ functions: see \cite{CRW} for this
characterization. Also,
this particular choice of $BMO$ function was exploited in \cite{JJ} in
their proof of weak convergence in $H^1$. It has
the following properties:

\begin{itemize}
\item $f \geq 0$ \item $\|f\|_{BMO} \leq \delta$ \item $f=1$ on
$E$
\end{itemize}

Observe that if $x \notin 2\Delta$, then
$M(h)(x) < \sigma(E)/\sigma(\Delta) < \eta$. For any
$\delta$, if we choose $\eta$ sufficiently small,
the function $1+\delta\log M(h)$ will be negative, and thus
$f=0$ outside $2\Delta$.

Using a standard mollification process (as in
\cite{Sa}) we can find a family $f_t$ of continuous
functions, $t>0$ such that:

\begin{itemize}
\item $f_t \rightarrow f$ in $L^p,$
\item For all t, there exists a $C$ such that $\|f_t\|_{BMO} \leq C\|f\|_{BMO},$
\item support of $f_t$ is contained in $3\Delta.$
\end{itemize}

\noindent Because
$f\ge 1$ on $E$, (\ref{eqn:normbound}) implies that

\begin{eqnarray}\nonumber
\frac{\omega(E)}{\omega(3\Delta)}&\le&
\omega(3\Delta)^{-1}\int_{3\Delta}
f\,d\omega=\omega(3\Delta)^{-1}\lim_{t\to0+}\int_{3\Delta}f_t\,d\omega\\
&\le& C_0\limsup_{t\to0+}\|f_t\|_{BMO}.\label{est1}
\end{eqnarray}

Hence by (\ref{est1})
$$\frac{\omega(E)}{\omega(3\Delta)}\le C_1\|f\|_{BMO}.$$

\bigskip
Now we choose $\delta$ so that $2C_1\delta < \varepsilon$, where
$C_1$ is the constant in the estimate above and this gives that

\begin{equation}
\frac{\omega(E)}{\omega(\Delta)} < M \epsilon
\end{equation}

\noindent where $M$ depends on the doubling constant of the
measure $\omega$.\vglue3mm

\medskip

Now that absolute continuity is established, the exact same argument gives
$A_{\infty}$. The function $f$, constructed in (\ref{eq:specialbmofunction}), will
have the same properties as before, except that, for general sets $E$, $f \geq 1$
a.e. $d\sigma$ on $E$, and hence a.e. $d\omega$ on $E$ by absolute continuity.

\bigskip

Before turning to the proof of the converse, we note the following corollary
of this argument.

\medskip

\noindent Suppose that the Dirichlet problem for $L$ with data in $L^p$ is
solvable in the sense that an apriori estimate in terms of square functions
holds:
$$\|S(u)\|_{L^p(\partial \Omega, d\sigma)} \lesssim \|f\|_{L^p(\partial \Omega, d\sigma)}.$$
Then the argument above shows that also
$$
\|N(u)\|_{L^p(\partial \Omega, d\sigma)} \lesssim
\|f\|_{L^p(\partial \Omega, d\sigma)}.$$

This can be derived from \ref{claim} as follows. Let $f$ be positive and supported in a surface
ball $\Delta$ of radius $r$, and let $\Delta'$ be as above. Then

\begin{equation}
(\omega(\Delta)^{-1}\int_{\Delta} f d\omega)^p \leq C
\sigma(\Delta)^{-1}\int_{\Delta'} S_r^p(u) d\sigma \leq C
\sigma(\Delta)^{-1}\int_{\Delta'} f^p d\sigma
\end{equation}

\noindent shows that $d\omega$
is absolutely continuous with respect to $d\sigma$ and the
density belongs to $B_q$, where $1/p+1/q = 1$.

\bigskip

\noindent{\it Proof of the Converse.}
This part of the proof of Theorem \ref{theorem:mainth} uses ideas in
in Fabes-Neri \cite{FN}, where the authors showed that
the $BMO$ Dirichlet problem was solvable for the
Laplacian in Lipschitz domains.

By
assumption, since $d\omega_L\in A_\infty(\delo,d\sigma)$,
there is $p_0>1$ such that the Dirichlet problem $(D_p)$
for $L$ is solvable for all $p_0<p\le\infty$.

Consider $f\in BMO(\delo)$. We will establish that
\begin{equation}
\iint_{T(\Delta)}|\nabla u|^2\delta(X)dX\le
C\sigma(\Delta)\|f\|^2_{BMO}.\label{0}
\end{equation}
 Consider any $\Delta\subset
\partial\Omega$ a surface ball of radius $r$. Let us denote by
$\widetilde\Delta$ and enlargement of $\Delta$ such that
$3\Delta\subset \widetilde\Delta\subset 5\Delta$. We will write
the solution $u$ of the Dirichlet problem for boundary data $f$ as
$u_1+u_2+u_3$, where $u_1, u_2$ solve
\begin{eqnarray}
Lu_1&=&0,\quad
u_1\big|_{\delo}=(f-f_{\widetilde\Delta})\chi_{\widetilde\Delta},\nonumber\\
Lu_2&=&0,\quad
u_2\big|_{\delo}=(f-f_{\widetilde\Delta})\chi_{\delo\setminus\widetilde\Delta},\nonumber\\
u_3&=&f_{\widetilde\Delta}\,\mbox{ in }\Omega.\nonumber
\end{eqnarray}
Here $f_{\widetilde\Delta}$ denotes, as before the average of $f$
over the set $\widetilde\Delta$ and $\chi_{\widetilde\Delta}$ is
the characteristic function of the set ${\widetilde\Delta}$.

We first estimate $u_1$. We claim that
\begin{equation}
\iint_{T(\Delta)}|\nabla u_1|^2\delta(X)dX\le
C\int_{{\widetilde\Delta}} S_r^2(u_1)\,d\sigma.\label{1}
\end{equation}

Let us denote by $\Delta_X$ the set $\{Q\in\partial\Omega;
X\in\Gamma(Q)\}$. It follows that $\sigma(\Delta_X\cap
\widetilde\Delta)\approx \delta(X)^{n-1}$. Hence
\begin{eqnarray}
\iint_{T(\Delta)}|\nabla u_1|^2\delta(X)dX&\le&
C\iint_{T(\Delta)}\delta(X)^{2-n}|\nabla u_1|^2\sigma(\Delta_X\cap
\widetilde\Delta)dX\nonumber \\&\le& C\int_{Q\in \widetilde\Delta}
\int_{\Gamma_r(Q)}\delta(X)^{2-n}|\nabla u_1|^2\,
dX\,d\sigma\label{2}\\
&\le& C\int_{{\widetilde\Delta}} S_r^2(u_1)\,d\sigma.\nonumber
\end{eqnarray}

By H\"older inequality for sufficiently large $p$ (such that the
$L^p$ Dirichlet problem is solvable on $\Omega$)
\begin{eqnarray}
\int_{{\widetilde\Delta}} S_r^2(u_1)\,d\sigma &\le&
\sigma(\widetilde\Delta)^{\frac{p-2}p}\left(\int_{{\widetilde\Delta}}
S^p(u_1)\,d\sigma\right)^{2/p}\nonumber\\
&\le&C\sigma(\Delta)^{\frac{p-2}p}\left(\int_{\widetilde\Delta}|u_1|^p\,d\sigma\right)^{2/p}.\label{3}
\end{eqnarray}
The last inequality uses solvability of the Dirichlet problem in
$L^p$, which implies that the $L^p$ norm the square function is
comparable to the $L^p$ norm of the boundary data. We put
(\ref{1}) and (\ref{3}) together to obtain an estimate
\begin{equation}
\iint_{T(\Delta)}|\nabla u_1|^2\delta(X)dX\le
C\sigma(\Delta)^{\frac{p-2}p}\left(\int_{\widetilde\Delta}|f-f_{\widetilde\Delta}|^p\,d\sigma\right)^{2/p}\le
C\sigma(\Delta)\|f\|^2_{BMO(p)}. \label{4}
\end{equation}
This is the desired estimate for $u_1$. Now we handle $u_2$. This
function is a solution of the equation $Lu_2=0$ with Dirichlet
boundary data
$f_2:=f-(f-f_{\widetilde\Delta})\chi_{\widetilde\Delta}$. Let us
call by $f_2^+$ and $f_2^-$ the positive and negative part of the
function $f_2$, that is $f_2=f_2^+-f_2^-$ and $f_2^+,f_2^-\ge 0$.
We denote by $u_2^{\pm}$ the solution of the Dirichlet problem
$$Lu_2^\pm=0\quad\mbox{in }\Omega,\quad u_2^\pm\big|_{\delo}=f_2^\pm.$$
Hence $u_2^\pm\ge 0$ and $u_2=u_2^+-u_2^-$. We claim the following

\begin{lemma}\label{l1} There exist $C>0$ depending only on the ellipticity
of the operator $L$ such that for any $X\in\Omega$
\begin{equation}
\left(\delta(X)^{-n}\int_{B(X,\,\delta(X)/2)}|\nabla
u_2^{\pm}(Y)|^2 dY\right)^{1/2}\le
\frac{C}{\delta(X)}\int_{\delo}f_2^{\pm}(Q)\,d\omega^X(Q).
\label{5}
\end{equation}
Here $\omega^X$ is the elliptic measure for the operator $L$ at
the point $X$.
\end{lemma}
This statement is a consequence of the Poincar\'e inequality that
allows to estimate the integral of a gradient by an average of
$(u_2^\pm-u_2^\pm(X))^2$ over slightly larger ball and by Harnack
inequality that implies $u_2^\pm(Y)\approx u_2^\pm(X)$ for $Y\in
B(X,\,\delta(X)/2)$. Notice that the integral $\int_{\delo}
f_2^{\pm}(Q)\,d\omega^X(Q)$ equals to the value of $u_2^\pm$ at
the point $X$.

Let us set
\begin{equation}
v_2(X)=\int_{\delo}|f_2(Q)|\,d\omega^X(Q)=\int_{\delo}(f_2^{+}(Q)+f_2^{-}(Q))\,d\omega^X(Q).
\label{6}
\end{equation}
It follows that $v^2(X)=u_2^+(X)+u_2^-(X)$.
\begin{lemma}\label{L-FN} There exist $C,\varepsilon>0$ depending only on the ellipticity constant of the operator $L$ such that for all $x\in T(\Delta)$:
\begin{itemize}
\item $v_2(X)\le C \|f\|_{BMO}$ \item $v_2(X)\le C
\|f\|_{BMO}\,\left(\frac{\delta(X)}r\right)^\varepsilon$. Here $r$
is the radius of the surface ball $\Delta$.
\end{itemize}
\end{lemma}
We postpone the proof of this lemma until we show how it gives us
the desired estimate.

To to that we consider a standard \lq dyadic' decomposition of the
Carleson region $T(\Delta)$. What this means is that $T(\Delta)$
can be written as a union of disjoint regions $I_n$,
$n=1,2,3,\dots$ such that for each region $I_n$ the diameter of
the region $d=\mbox{diam}(I_n)$ is comparable to the distance
$\mbox{dist}(I_n,\delo)$ and the volume of the region is
comparable to $d^n$. For each region $I_n$ we denote by $x_n$ a
point inside $I_n$. It follows that
\begin{eqnarray}
&&\iint_{T(\Delta)}|\nabla u_2^\pm|^2\delta(X)dX\le \sum_n
\iint_{I_n}|\nabla u_2^\pm|^2\delta(X)dX\le C\sum_n
\delta(x_n)\frac{|u_2^\pm(x_n)|}{\delta(x_n)^2}\delta(x_n)^n\nonumber\\
&\le& C\int_{T(\Delta)}\frac{|u_2^{\pm}(X)|^2}{\delta(X)}dX\le
C\left(r^{-2\varepsilon}\int_{T(\Delta)}\delta(X)^{2\varepsilon-1}dX\right)\|f\|^2_{BMO}.\label{7}
\end{eqnarray}
Here we used Lemma \ref{l1} for the last estimate in the first
line of (\ref{7}) and Lemma \ref{L-FN} for the last estimate in
the second line (clearly $u_2^\pm(X)\le v_2(X)$).

Since
$r^{-2\varepsilon}\int_{T(\Delta)}\delta(X)^{2\varepsilon-1}dX\le
Cr^{n-1}\approx\sigma(\Delta)$ we see that (\ref{4}) and (\ref{7})
together implies the estimate (\ref{0}) we sought (function $u_3$
is constant, hence the required estimate hold trivially).
\vglue2mm

\noindent{\it Proof of Lemma \ref{L-FN}.} The first estimate of
the lemma, namely that $v_2(X)\le C\|f\|_{BMO}$, essentially follows from Lemma
on p.35 in \cite{FN}. As stated there
$$v_2(X)=\int_{\delo\setminus \widetilde\Delta}|f-f_{\widetilde\Delta}|K(X,Q)d\sigma(Q),$$
for some kernel $K(X,Q)$ (a Radon-Nykodim derivative of the
elliptic measure $\omega^X$). Fabes and Neri then use then fact
that $K\in B^2(d\sigma)$\footnote{We denote by $B^q$ the class of
Gehring weights. The weights in this class satisfy the reverse
H\"older inequality with exponent $q$.} to establish the estimate.
By looking at their proof we see that it is enough to have $K\in
B^q$ for some $q>1$. This holds, as we assume that $\omega^X\in
A_\infty(d\sigma)=\bigcup_{q>1} B_q(d\sigma)$.

The further improvement in the estimate $v_2(X)\le
C\|f\|_{BMO}\,\left(\frac{\delta(X)}r\right)^\varepsilon$ is a
consequence of Di Giorgi-Nash-Moser theory. Nonnegative solutions
$u$ of $L$ in the region $T(\widetilde\Delta)$ which vanish on
$2\Delta$ satisfy
$$u(X)\le C\left(\frac{|X-Q|}{r}\right)^\varepsilon \sup_{T(2\Delta)} u,\qquad\mbox{for any }X\in T(\Delta).$$
Here $\varepsilon$ only depends on the ellipticity constant of the
operator $L$ and $Q$ is the center of the ball $\Delta$. (See for
example (1.9) in \cite{KKPT} for reference). From this the
estimate follows as we can move point $Q$ around (within $\Delta$)
as our function vanishes on $\widetilde\Delta\supset 3\Delta$.
\vglue2mm

Now we prove the reverse estimate to (\ref{0}). We want to show
that
\begin{equation}
\|f\|^2_{BMO^*(d\sigma)}\le C\sup_{\Delta\subset
\delo}\iint_{T(\Delta)}|\nabla
u|^2\delta(X)\frac{dX}{\sigma(\Delta)}.\label{9}
\end{equation}
In this case it is more convenient to use (\ref{8}) to define BMO
norm. We first prove the following
\begin{equation}
\sup_{\Delta\subset\delo}\inf_{c_{\Delta}}
\sigma(\Delta)^{-1}\int_{\Delta}|f-c_\Delta|d\sigma \le C
\sup_{\Delta\subset\delo}\left(\inf_{c_{\Delta}}
\omega(\Delta)^{-1}\int_{\Delta}|f-c_\Delta|^p
d\omega\right)^{1/p}. \label{10}
\end{equation}
Here $\omega=\omega^{X_0}$ is the elliptic measure for the
operator $L$ at some (fixed) interior point $X_0$.
 This inequality implies that a BMO function with respect to the
surface measure $\sigma$ is also a BMO function with respect to
the elliptic measure $\omega$. Indeed, Let $d\sigma=kd\omega$. The
fact $\omega\in A_\infty(d\sigma)$ implies that $\sigma\in
A_\infty(d\omega)=\bigcup_{q>1}B_q(d\omega)$. Hence there exists
$q>1$ such that $k$ satisfies the reverse H\"older inequality
\begin{equation}
\left(\omega(\Delta)^{-1}\int_{\Delta}k^q\,
d\omega\right)^{1/q}\le C\omega(\Delta)^{-1}\int_{\Delta}k\,
d\omega\quad\mbox{for all }\Delta\subset\delo. \label{11}
\end{equation}
It follows
\begin{eqnarray}
&&\sigma(\Delta)^{-1}\int_{\Delta}|f-c_\Delta|d\sigma
=\sigma(\Delta)^{-1}\int_{\Delta}|f-c_\Delta|k\,d\omega\nonumber\\&\le&
\sigma(\Delta)^{-1}\left(\int_{\Delta}k^q\,d\omega\right)^{1/q}
\left(\int_{\Delta}|f-c_\Delta|^pd\omega\right)^{1/p}\label{12}\\
&\le&C\sigma(\Delta)^{-1}\omega(\Delta)^{1/q-1}\left(\int_{\Delta}k\,d\omega\right)
\left(\int_{\Delta}|f-c_\Delta|^pd\omega\right)^{1/p}.\nonumber\\
&=&C
\left(\omega(\Delta)^{-1}\int_{\Delta}|f-c_\Delta|^pd\omega\right)^{1/p}.
\end{eqnarray}
This gives (\ref{10}). It also follows that it suffices to prove
(\ref{9}) with $d\omega$ measure on the left-hand side instead of
$d\sigma$.

In what follows we use the following lemma from \cite{K}.
\begin{lemma}\label{l2} Let $X_0$ be a fixed point inside a
Lipschitz domain $\Omega$, $\omega^{X_0}$ the elliptic measure for
an operator $L$ at $X_0$ and $G(.,.)$ the Green's function for
$L$. Then for any open surface ball $\Delta_r\subset \delo$ or
radius $r$ such that $\delta(X_0)\ge 2r$ and
\begin{equation}
G(X_0,Y)r^{n-2}\approx \omega(\Delta_r),\label{13}
\end{equation}
where $Y\in\Omega$ such that dist$(Y,\Delta_r)\approx
\delta(Y)=r$. The precise constants in the estimate (\ref{13})
only depends on the ellipticity of $L$ and Lipschitz character of
domain $\Omega$.
\end{lemma}

The following lemma is crucial for the proof.

\begin{lemma}\label{l3}
There exists $C>0$ such that for all $f\in BMO(d\omega)$
\begin{equation}
 \|f\|_{BMO^*(d\omega)}\le C\sup_{\Delta\subset \delo}\left(\iint_{T(\Delta)}|\nabla u|^2G(X_0,X)\frac{dX}{\omega(\Delta)}\right)^{1/2}.\label{14}
\end{equation}
\end{lemma}
Assume for the moment the Lemma is true. By using Lemma \ref{l2}
we get that
\begin{equation}
\iint_{T(\Delta)}|\nabla u|^2G(X_0,X)dX\le
C\iint_{T(\Delta)}|\nabla u|^2\delta(X)^{2-n}\omega(\Delta_X) dX,
\end{equation}
where $\Delta_X$ is as before the set $\{Q\in\delo;
X\in\Gamma(Q)\}$. By changing the order of integration we get that
\begin{equation}
\iint_{T(\Delta)}|\nabla u|^2\delta(X)^{2-n}\omega(\Delta_X) dX\le
\int_{\widetilde\Delta} S^2_ru(Q)\,d\omega(Q).\label{15}
\end{equation}
Combining (\ref{14})-(\ref{15}) we get that
\begin{equation}
 \|f\|_{BMO^*(d\omega)}\le \sup_{\Delta\subset\delo}\left(\int_{\Delta} S^2_ru(Q)\,\frac{d\omega(Q)}{\omega(\Delta)}\right)^{1/2}. \label{16}
\end{equation}
Now we use the same trick as above to change measure back from
$\omega$ to $\sigma$. Again using reverse H\"older inequality (now
for $k^{-1}$) we get that
\begin{equation}
\sup_{\Delta_r\subset\delo}\left(\int_{\Delta_r}
S^2_ru(Q)\,\frac{d\omega(Q)}{\omega(\Delta_r)}\right)^{1/2}\le
C\sup_{\Delta_r\subset\delo}\left(\int_{\Delta_r}
S^q_ru(Q)\,\frac{d\sigma(Q)}{\sigma(\Delta_r)}\right)^{1/q}\quad\mbox{for
some }q>2.\nonumber
\end{equation}
Finally, there  exists $C>0$
\begin{eqnarray}
\sup_{\Delta_r\subset\delo}\left(\int_{\Delta_r}
S^q_ru(Q)\,\frac{d\sigma(Q)}{\sigma(\Delta_r)}\right)^{1/q}&\le&
C\sup_{\Delta_r\subset\delo}\left(\int_{\Delta_r}\label{17}
S^2_ru(Q)\,\frac{d\sigma(Q)}{\sigma(\Delta_r)}\right)^{1/2}\\&=&C\sup_{\Delta\subset
\delo}\left(\iint_{T(\Delta)}|\nabla
u|^2\delta(X)\frac{dX}{\sigma(\Delta)}\right)^{1/2}.\nonumber
\end{eqnarray}
The first estimate in (\ref{17}) follows from the BMO
John-Nirenberg argument (same way as (\ref{88}) is established).
This concludes the proof of Theorem \ref{theorem:mainth} (modulo
Lemma \ref{l3}).\vglue2mm

\noindent {\it Proof of Lemma \ref{l3}.} We fix a surface ball
$\Delta\subset \delo$ or radius $r$ and center $Q$. As before we
consider a point $X_0$ inside $\Omega$ such that $\delta(X_0)\ge
5r$. Finally, let us denote by $\cal D$ the domain $\Omega \cap
B(Q,4r)$. We pick a point $X\in \cal D$ such that
dist$(X,\partial{\cal D})\approx 2r$. We denote by $\nu$ the
elliptic measure for operator $L$ on the domain $\cal D$ with pole
at $X$.

We study relations between measures $\omega$ and $\nu$. The
following Lemma holds

\begin{lemma}\label{l4} For any measurable set $E\subset\Delta$
\begin{equation}
\label{18}\frac{\omega(E)}{\omega(\Delta)}\le C\nu(E),
\end{equation}
where the constant $C>0$ only depends on the ellipticity constant
and Lipschitz character of the domain $\Omega$.
\end{lemma}

It suffices to establish (\ref{18}) for all balls
$\Delta'\subset\Delta$, as the general statement for all
measurable sets $E$ follows by a covering lemma. For both balls
$\Delta'$ and $\Delta$ we find points $Y'$ and $Y$, respectively
such that dist$(Y',\partial\Delta')\approx\delta(Y')=r'$ and
dist$(Y,\partial\Delta)\approx\delta(Y)=r$, where $r'$ and $r$ are
radii of these balls. According to Lemma \ref{l2}
$$\omega(\Delta')\approx G_{\Omega}(X_0,Y')(r')^{n-2},\quad\mbox{and}\quad \nu(\Delta')\approx G_{\cal D}(X,Y')(r')^{n-2}.$$
Hence
$$\frac{\omega(\Delta')}{\nu(\Delta')}\approx\frac{G_{\Omega}(X_0,Y')}{G_{\cal D}(X,Y')}\approx \frac{G_{\Omega}(X_0,Y)}{G_{\cal D}(X,Y)}.$$
The last relation comes form the comparison principle for two
positive solutions $v(.)=G_{\Omega}(X_0,.)$ and $w(.)=G_{\cal
D}(X,.)$ that vanish at the boundary. Finally,
$$\frac{\omega(\Delta')}{\nu(\Delta')}\approx\frac{G_{\Omega}(X_0,Y)}{G_{\cal D}(X,Y)}\approx\frac{\omega(\Delta)r^{n-2}}{\nu(\Delta)r^{n-2}},$$
again by using Lemma \ref{l2}. However, $\nu(\Delta)=O(1)$, since
the measure $\nu$ is doubling, and $\nu(\partial{\cal D})=1$.
Hence Lemma \ref{l4} follows.\vglue2mm

By Lemma \ref{l4} we see that for any $c_\Delta\in \R$
\begin{equation}
\label{19}\int_{\Delta}|f-c_{\Delta}|^2\frac{d\omega}{\omega(\Delta)}\le
C \int_{\Delta}|f-c_{\Delta}|^2 d\nu\le C\int_{\partial{\cal
D}}|u-c_\Delta|^2 d\nu.
\end{equation}
Since $\nu$ is the natural (elliptic) measure for the domain $\cal
D$ it follows that the $L^2(d\nu)$ Dirichlet problem is always
solvable in this domain. This implies the the $L^2(d\nu)$ norm of
the square function is comparable with the $L^2(d\nu)$ of the
(normalized) boundary data, i.e.,
\begin{equation}
\label{20}\inf_{c_{\Delta}\in\R}\int_{\partial{\cal
D}}|u-c_\Delta|^2 d\nu\approx \int_{\partial{\cal D}} S^2u
\,d\nu\approx \iint_{\Omega\setminus B_{r/8}(X)}|\nabla
u(Y)|^2G_{\cal D}(X,Y)dY.
\end{equation}
Finally, we claim that
\begin{equation}
\label{21}G_{\cal D}(X,Y)\le G_{\Omega}(X,Y)\approx
\frac{G_{\Omega}(X_0,Y)}{\omega(\Delta)},\quad\mbox{for all }Y\in
\Omega\setminus B_{r/8}(X).
\end{equation}
Combining the estimates (\ref{19})-(\ref{21}) we obtain Lemma
\ref{l3}. The first estimate of (\ref{21}) is simply a maximum
principle, as $G_{\cal D}(X,Y)$ vanishes on the whole
$\partial\cal D$, and $G_{\Omega}(X,Y)$ is positive at the portion
of this boundary. Both functions have same pole at $X$. The
relation $G_{\Omega}(X,Y)\approx
\frac{G_{\Omega}(X_0,Y)}{\omega(\Delta)}$ can be established as
follows. For $Y\in \Omega\setminus B_{r/8}(X)$ such that
$\delta(Y)\ge r$ Lemma \ref{l2} implies that
$G_{\Omega}(X_0,Y)\approx r^{n-2}\omega(\Delta)$. On the other
hand $G_{\Omega}(X,Y)\approx r^{n-2}$ as $Y$ is of distance $r$
from the pole and also $r$ away from the boundary. For $Y$ near
the boundary we use the comparison principle (since both function
vanish at $\delo$. This gives
$$\frac{G_{\Omega}(X,Y)}{G_{\Omega}(X_0,Y)}\approx \frac{G_{\Omega}(X,Y')}{G_{\Omega}(X_0,Y')}$$
for all $Y, Y'\in \Omega\setminus B_{r/8}(X)$. This establishes
(\ref{21}) and concludes the proof of Theorem
\ref{theorem:mainth}.\qed \vglue2mm

\noindent{\it Proof of Theorem\ref{theorem:mainth2}.} By Theorem
\ref{theorem:mainth} if follows that $d\omega_L\in
A_\infty(\partial\Omega,d\sigma)$. Since
$$A_\infty(\partial\Omega,d\sigma)=\bigcup_{p>1}B_p(\partial\Omega,d\sigma),$$
we see that $d\omega_L\in B_p(\partial\Omega,d\sigma)$ for some
$p>1$.  From this the claim follows since $d\omega_L\in
B_p(\partial\Omega,d\sigma)$ implies the solvability of the
$L^{p'}$ Dirichlet problem. The range of solvability
$(p_0,\infty)$ can be then obtained by realizing that
$B_p(\partial\Omega,d\sigma)\subset B_q(\partial\Omega,d\sigma)$
for $q<p$.\qed \vglue 2mm

\noindent{\it Proof of the Remark 2.1.} Indeed, by Theorem
\ref{theorem:mainth2} given that (\ref{equiv}) holds, the $L^p$
Dirichlet boundary value problem is solvable for some large
$p<\infty$. Consider now an arbitrary BMO function
$f:\partial\Omega\to\R$. As we argue in (\ref{est1}), there exists
a sequence of {\it continuous} functions $f_n:\partial\Omega\to\R$
such that $f_n\to f$ in $L^p(\partial\Omega)$ and
$\|f_n\|_{BMO}\le C\|f\|_{BMO}$ for some $C>0$ independent of $n$.

For each $f_n$ we can solve the continuous Dirichlet boundary
value problem which will give us solutions $u_n$ such that
$$\|N(u_n)\|_{L^p(\partial\Omega)}\le C\|f_n\|_{L^p(\partial\Omega)}\le C\|f\|_{BMO}.$$
In addition, also
$$\|N(u_n-u_m)\|_{L^p(\partial\Omega)}\le
C\|f_n-f_m\|_{L^p(\partial\Omega)}\to 0,\qquad\mbox{as
}n,m\to\infty,$$ since $f_n\to f$ in $L^p$ and (\ref{dp}) holds.
This implies that the sequence $(u_n)_{n\in\N}$ is locally
uniformly Cauchy in $L^\infty_{loc}(\Omega)$, hence
$$u(X)=\lim_{n\to\infty} u_n(X),\qquad\mbox{for }X\in\Omega$$
is pointwise well defined.

We claim that this $u$ is a weak solution to $Lu=0$. That is,

\begin{equation}\label{d-sol2}
\int_{\Omega}A(X)\nabla u(X).\nabla\psi(X)\,dX=0,\qquad\mbox{for
all }\psi\in C_0^\infty(\Omega),
\end{equation}

To see this, fix a compact set $K\subset\Omega$. By the dominated
convergence theorem we know that
$$u_n\to u,\qquad\mbox{ in any }L^p(K),\quad p<\infty.$$
Hence for any $K'\subset\subset K$ by Cacciopoli we have that
$$\int_{K'}|\nabla(u_n-u_m)(X)|^2\,dX\le C_{K,K'}\int_{K}|(u_n-u_m)(X)|^2\,dX\to 0,\quad\mbox{as }n,m\to\infty.$$
It follows that $\nabla u_n$ converges locally uniformly in $L^2$,
from which we get that $u$ belongs to $W^{1,2}_{loc}(\Omega)$ and
$\nabla u_n\to \nabla u$ in $L^2_{loc}(\Omega)$. Therefore
(\ref{d-sol2}) follows as we already now that (\ref{d-sol2}) holds
for every $u_n$ and we can pass to the limit $n\to\infty$.

Hence with the use of Fatou's lemma (see Appendix B of \cite{Di}
for details) we get that $N(u-u_n)\to 0$ in $L^p(\partial\Omega)$
as $n\to\infty$. This implies that $\|N(u)\|_{L^p}<\infty$, so
$N(u)(Q)<\infty$ a.e. for $Q\in\partial\Omega$ and also one has
existence of nontangential limits a.e.: $\lim_{X\to Q,\,
X\in\Gamma(Q)} u(X)$.

Finally, we also get that (\ref{equiv}) will also hold for $u$ by
the limiting argument, since it holds for each $u_n$:
\begin{equation}
\sup_{\Delta\subset\delo}\sigma(\Delta)^{-1}\iint_{T(\Delta)}|\nabla
u_n|^2\delta(X)dX \lesssim \sup_{I\subset
\delo}\sigma(I)^{-1}\int_I |f_n-f_{n,I}|^2 d\sigma. \label{equiv2}
\end{equation}
Notice that taking the limsup on the right-hand side of
(\ref{equiv2}) yields just a multiple of BMO norm of $f$, as
$\|f_n\|_{BMO}\le C\|f\|_{BMO}$. On the left-hand side we may take
the limit $$\sigma(\Delta)^{-1}\iint_{T(\Delta)\setminus {\cal
C}_\varepsilon}|\nabla u_n|^2\delta(X)dX\to
\sigma(\Delta)^{-1}\iint_{T(\Delta)\setminus {\cal
C}_\varepsilon}|\nabla u|^2\delta(X)dX, \quad n\to\infty,$$ since
$\nabla u_n\to\nabla u$ in $L^2_{loc}(\Omega)$. Here ${\cal
C}_{\varepsilon}=\{X\in\Omega:\,
\mbox{dist}(X,\partial\Omega)<\varepsilon\}$. It follows that for
any $\varepsilon>0$
\begin{equation}
\sup_{\Delta\subset\delo}\sigma(\Delta)^{-1}\iint_{T(\Delta)\setminus
{\cal C}_\varepsilon}|\nabla u|^2\delta(X)dX \lesssim
\|f\|^2_{BMO}. \label{equiv3}
\end{equation}
As the constant in (\ref{equiv3}) does not depend on $\varepsilon$
we get the required estimate on the whole $T(\Delta)$. In fact, it
can be shown that \emph {equivalence} holds between the two
quantities in (\ref{equiv}).\qed

%%%%%%%%%%%%%%%%%%%%%%%%%%%%%%%%%%%%%%%%%%%%%%%%%%%%%%%%%%%%%%%%%%%%%%%%%%%%%%%%%%%%%%%%%%%%%%%%%%%%%%%%

\end{document}